\documentclass[12pt]{article}
\usepackage[dvipsnames]{xcolor}
\usepackage{amsmath}
\usepackage{setspace}
\usepackage{bm}
\usepackage{bbm}
\usepackage{amssymb}
\usepackage{indentfirst}
\usepackage[superscript]{cite}
\usepackage{geometry}
\usepackage{mathtools}

\newcommand*{\citena}[1]{%
\begingroup
[\color{Green}
\romannumeral-`\x 
\setcitestyle{numbers}%
\cite{#1}%
\endgroup
]\ignorespacesafterend
}

\newcommand*{\citesup}[1]{%
\begingroup
\color{Green}
\cite{#1}%
\endgroup
\ignorespacesafterend
}

\newcommand*{\eqrefe}[1]{%
\begingroup
(\color{BrickRed}
\romannumeral-`\x 
\setcitestyle{numbers}%
\ref{eq:#1}%
\endgroup
)\ignorespacesafterend
}

\newcommand*{\secrefe}[1]{%
\begingroup
(\color{Aquamarine}
\romannumeral-`\x 
\setcitestyle{numbers}%
\ref{#1}%
\endgroup
)\ignorespacesafterend
}

\newtagform{Tags}[\textcolor{BrickRed}]{\color{Black}(}{)}

\usepackage[super]{natbib}
\setcitestyle{super}

\newcommand{\ii}{\bm{i}}

\DeclareMathOperator{\csch}{csch}
\DeclareMathOperator{\HP}{HP}

\begin{document}
\title{The Hurwitz Zeta Function at the Positive Integers}
\author{Jose Risomar Sousa}
\date{February 19, 2019}
\maketitle
\usetagform{Tags}

\begin{abstract}
A formula for the Hurwitz zeta function at the positive integers $k$, $\zeta(k,b)$, is created by solving the real and the imaginary parts separately and then combining them. A few different formulae for the Hurwitz zeta function are known from the literature, but they are very general and usually hold for $\Re{(k)}>1$. The advantage of formulae that only hold at the positive integers is the fact that they are simpler and easier to work with. An analytic continuation of the generating function of $\zeta(k,b)$ is also obtained as $\sum_{k\ge 2}x^k(\zeta(k,b)-1/b^k)$, where the term $1/b^k$ was subtracted for convenience.
\end{abstract}

\tableofcontents

\section{Introduction}
The Hurwitz zeta function, $\zeta(k,b)$, is one of the many generalizations that were conceived for the Riemann zeta function. Its importance lies on its relation to the Riemann zeta function and consequently to the Riemann hypothesis itself to some degree.\\

In this manuscript, a formula for $\zeta(k,b)$ is created that holds at the positive integers $k$ greater than one. The advantage of formulae that only hold at the positive integers is the fact that they are simpler and easier to work with. It is an obvious statement if, for example, one compares the closed form of the zeta function at the positive even integers with its general integral, valid for $\Re{(k)}>1$.\\

Findings from previous papers released by the author on generalized harmonic numbers, $H_k(n)$\citesup{GHNR}, and generalized harmonic progressions, $\HP_k(n)$\citesup{GHP}$^,$\citesup{CHP}, make it relatively easy to determine the limit of $\HP_k(n)$ as $n$ goes to infinity, which is $\zeta(k,b)$.\\

For the problem under consideration, one must rely on the formula for $\HP_k(n)$ from \citena{CHP}, since it allows for non-integer $b$.\\

First, let us recall the formula, which holds for any complex $b$, except for a zero measure subset of $\mathbb{C}$ (that is, $\ii\,b \in \mathbb{Z}$). For simplicity, let us assume that $a=1$ and that $b$ is real, so as to know which part of the formula is real and which part is purely complex. Then, $\HP_k(n)$ is given by:
\begin{multline} \nonumber
\sum_{j=1}^{n}\frac{1}{(\ii j+b)^k} =-\frac{1}{2b^k}+\frac{1}{2(\ii n+b)^k}\\ +(2\pi)^k e^{-2\pi b}\int_{0}^{1}\left(\frac{(1-u)^{k-1}}{(k-1)!}+\sum_{j=1}^{k}\frac{\mathrm{Li}_{-j+1}(e^{-2\pi b})(1-u)^{k-j}}{(j-1)!(k-j)!}\right)e^{\pi u(\ii n+2b)}\sin{\pi n u}\cot{\pi u}\,du \text{,}
\end{multline}
\noindent where $\mathrm{Li}_{-j+1}(e^{-2\pi b})$ is the polylog of order $-j+1$. The polylog is the analytic continuation of the Dirichlet series below:
\begin{equation} \nonumber
\mathrm{Li}_{k}(z)=\sum_{j=1}^{\infty}\frac{z^j}{j^k}
\end{equation}
\indent The polynomial in $u$ that goes within the integral of $\HP_k(n)$ is generated by the below function:
\begin{equation} \nonumber
f(x)=-\frac{2\pi x\,e^{2\pi x(1-u)}}{e^{2\pi x}-e^{2\pi b}} \Rightarrow \frac{f^{(k)}(0)}{k!}=(2\pi)^k e^{-2\pi b}\left(\frac{(1-u)^{k-1}}{(k-1)!}+\sum_{j=1}^{k}\frac{\mathrm{Li}_{-j+1}(e^{-2\pi b})(1-u)^{k-j}}{(j-1)!(k-j)!}\right)
\end{equation}
\indent Assuming $b$ is real, this integral can be transformed using Euler's formula for the exponential of a complex argument, and then replacing $\cos{\pi nu}\sin{\pi nu}$ and $\sin^2{\pi nu}$ with equivalent expressions. The Kronecker delta ($\delta_{ij}=1$ if $i=j$, and 0 otherwise) is also introduced into the formula, to make it shorter:
\begin{multline} \nonumber
\sum_{j=1}^{n}\frac{1}{(\ii j+b)^k} =-\frac{1}{2b^k}+\frac{1}{2(\ii n+b)^k}\\+\frac{(2\pi)^k}{2}\int_{0}^{1}\sum_{j=1}^{k}\frac{\left(\delta_{1j}+\mathrm{Li}_{-j+1}(e^{-2\pi b})\right)u^{k-j}}{(j-1)!(k-j)!}e^{-2\pi b u}\left(\sin{2\pi n(1-u)}+\ii(1-\cos{2\pi n (1-u)})\right)\cot{\pi(1-u)}\,du \text{}
\end{multline} 

\section{The limit of $\HP_k(n)$ for real $b$} 
In this section the real and imaginary parts of the limit of $\HP_k(n)$ as $n$ approaches infinity are obtained separately, assuming $b$ is real. At the end, $\zeta(k,-\ii\,b)$ is obtained by means of the relation:
\begin{equation} \label{eq:relacao}
\zeta(k,-\ii\,b)=\ii^k\sum_{j=0}^{\infty}\frac{1}{(\ii j+b)^k} \text{}
\end{equation} 
\indent This solution provides another proof that $\HP(n)$ diverges, though that is not the main focus herein.

\subsection{The real part of $\zeta(k,-\ii\,b)$} \label{Real}
For the real part, when $n$ is large one has:
\begin{multline} \nonumber
\Re\left(\sum_{j=1}^{n}\frac{1}{(\ii j+b)^k}\right) \\\sim -\frac{1}{2b^k}+\frac{(2\pi)^k}{2}\int_{0}^{1}\sum_{j=1}^{k}\frac{\left(\delta_{1j}+\mathrm{Li}_{-j+1}(e^{-2\pi b})\right)u^{k-j}}{(j-1)!(k-j)!}e^{-2\pi b u}\sin{2\pi n(1-u)}\cot{\pi(1-u)}\,du 
\end{multline} 
\indent Let us recall a result that appeared in reference \citena{GHNR}, whose proof depends on formulae that feature in Abramowitz and Stegun\citesup{Abramo}:
\begin{equation}
\textbf{Theorem 1}\lim_{n\to\infty}\int_{0}^{1}u^{k}\sin{2\pi n(1-u)}\cot{\pi(1-u)}\,du=
\begin{cases} \nonumber
      1, & \text{if}\ k=0\\
      1/2, & \text{if integer }k>0
\end{cases} 
\end{equation}
\indent In the formula of the real part of $\zeta(k,-\ii\,b)$, the summation within the integral can be split into a polynomial in $u$ (without intercept) and a constant, to handle the fact that the limits for each are different. And since $e^{-2\pi bu}$ is a polynomial in $u$ itself, when they are multiplied together, as below, in only one of the four possible cases the limit as $n$ approaches infinity is $1$ (all others are $1/2$):
\begin{multline} \nonumber
\lim_{n\to\infty}\int_{0}^{1}\left(\frac{\delta_{1k}+\mathrm{Li}_{-k+1}(e^{-2\pi b})}{(k-1)!}+\sum_{j=1}^{k-1}\frac{\left(\delta_{1j}+\mathrm{Li}_{-j+1}(e^{-2\pi b})\right)u^{k-j}}{(j-1)!(k-j)!}\right)\\ \cdot \left(1+\sum_{q=1}^{\infty}\frac{(-2\pi b)^q u^q}{q!}\right)\sin{2\pi n(1-u)}\cot{\pi(1-u)}\,du 
\end{multline}\\
\indent This yields the below limit:
\begin{equation} \nonumber
\frac{\delta_{1k}+\mathrm{Li}_{-k+1}(e^{-2\pi b})}{(k-1)!}\left(1+\frac{-1+e^{-2\pi b}}{2}\right)+\frac{e^{-2\pi b}}{2}\sum_{j=1}^{k-1}\frac{\left(\delta_{1j}+\mathrm{Li}_{-j+1}(e^{-2\pi b})\right)u^{k-j}}{(j-1)!(k-j)!}
\end{equation}
\indent Therefore, as a consequence of Theorem 1, after all the necessary calculations are carried out, one concludes that if $b$ is real then:
\begin{multline} \label{eq:Hurwitz_real_part}
\Re{\left(\sum_{j=0}^{\infty}\frac{1}{(\ii j+b)^k}\right)}\\=\frac{1}{2b^k}+\frac{(2\pi)^k\left(\delta_{1k}+\mathrm{Li}_{-k+1}(e^{-2\pi b})\right)}{4(k-1)!}+\frac{(2\pi)^ke^{-2\pi b}}{4}\sum_{j=1}^{k}\frac{\delta_{1j}+\mathrm{Li}_{-j+1}(e^{-2\pi b})}{(j-1)!(k-j)!}
\end{multline}
\indent Note the summation index $j$ starts at zero to coincide with the Hurwitz zeta function, so it is necessary to add $1/b^k$.\\

The real part of this infinite series is finite even when $k=1$, unlike the imaginary part. Next, the same reasoning is applied for the imaginary part.

\subsection{The imaginary part of $\zeta(k,-\ii\,b)$} \label{approx}
For the imaginary part, when $n$ is large one has:
\begin{multline} \nonumber \label{eq:im_part_1}
\Im\left(\sum_{j=1}^{n}\frac{1}{(\ii j+b)^k}\right) \\ \sim \frac{(2\pi)^k}{2}\int_{0}^{1}\sum_{j=1}^{k}\frac{\left(\delta_{1j}+\mathrm{Li}_{-j+1}(e^{-2\pi b})\right)u^{k-j}}{(j-1)!(k-j)!}e^{-2\pi b u}(1-\cos{2\pi n (1-u)})\cot{\pi(1-u)}\,du 
\end{multline} 
\indent To fully understand what happens in the case of the imaginary part, one needs to go back to one of the results from \citena{GHNR}. It provides an equivalence between certain integrals and generalized harmonic numbers. More specifically, it has been proved that,
\begin{multline} \nonumber
\int_{0}^{1}(1-u)^{2k+1}\left(1-\cos{2\pi n u}\right)\cot{\pi u}\,du \\ =\frac{2(-1)^k(2k+1)!}{(2\pi)^{2k+1}}\left(\sum_{j=0}^{k}\frac{(-1)^{k-j}(2\pi)^{2k-2j}}{(2k+1-2j)!}H_{2j+1}(n)-\frac{1}{2 n^{2k+1}}\sum_{j=0}^{k}\frac{(-1)^{j}(2\pi n)^{2j}}{(2j+1)!}\right) \text{,}
\end{multline}
\noindent and
\begin{multline} \nonumber
\int_{0}^{1}(1-u)^{2k+2}\left(1-\cos{2\pi n u}\right)\cot{\pi u}\,du \\ =\frac{2(-1)^k(2k+2)!}{(2\pi)^{2k+1}}\left(\sum_{j=0}^{k}\frac{(-1)^{k-j}(2\pi)^{2k-2j}}{(2k+2-2j)!}H_{2j+1}(n)-\frac{1}{2 n^{2k+1}}\sum_{j=0}^{k}\frac{(-1)^{j}(2\pi n)^{2j}}{(2j+2)!}\right) 
\end{multline}\\
\indent Let us consider their limit as $n$ approaches infinity. If $n$ is sufficiently large:
\begin{equation} \nonumber
\int_{0}^{1}(1-u)^{2k+1}\left(1-\cos{2\pi n u}\right)\cot{\pi u}\,du \sim \frac{\gamma+\log{n}}{\pi}+\frac{2(-1)^k(2k+1)!}{(2\pi)^{2k+1}}\sum_{j=1}^{k}\frac{(-1)^{k-j}(2\pi)^{2k-2j}}{(2k+1-2j)!}\zeta{(2j+1)} 
\end{equation}
\begin{equation} \nonumber
\int_{0}^{1}(1-u)^{2k+2}\left(1-\cos{2\pi n u}\right)\cot{\pi u}\,du \sim \frac{\gamma+\log{n}}{\pi}+\frac{2(-1)^k(2k+2)!}{(2\pi)^{2k+1}}\sum_{j=1}^{k}\frac{(-1)^{k-j}(2\pi)^{2k-2j}}{(2k+2-2j)!}\zeta(2j+1) \text{}
\end{equation}
\indent Using findings from \citena{GHNR}, the following analytic continuation can be obtained for these approximations\footnote{It stems from the expressions for $C^z_{2k+1}(n)$ and $S^z_{2k}(n)$ and for their limits, $C^z_{2k+1}$ and $S^z_{2k}$.}, which must hold for complex $k$ such that $\Re{(k)}>0$:
\begin{equation} \nonumber
\int_{0}^{1}(1-u)^k(1-\cos{2\pi n u})\cot{\pi u}\,du \sim \frac{\gamma+\log{n}}{\pi}-\int_{0}^{1}(u^k-u)\cot{\pi u}\,du
\end{equation}
\indent One can see that, unlike the integrals that appear in the formula of the real part, these new ones all diverge as $n$ goes to infinity.\\ 

Therefore, like in Theorem 7 from \citena{GHNR}, a linear combination of these integrals, say $p(u)=\sum_{k}a_k u^k$,
\begin{equation} \nonumber
\int_{0}^{1}p(u)(1-\cos{2\pi n(1-u)})\cot{\pi(1-u)}\,du \text{,}
\end{equation}\\
\noindent only converges if $p(1)=0$. But that is not a problem, since fortunately,
\begin{equation} \nonumber
\int_{0}^{1}(1-\cos{2\pi n(1-u)})\cot{\pi(1-u)}\,du=0 \text{, } \forall \text{ integer }n\text{,}
\end{equation}
\noindent which means that, without changing the result whatsoever for integer $n$, the formula at the beginning of this section can be changed to:
\begin{equation} \nonumber
\frac{(2\pi)^k}{2}\int_{0}^{1}\sum_{j=1}^{k}\frac{\left(\delta_{1j}+\mathrm{Li}_{-j+1}(e^{-2\pi b})\right)\left(u^{k-j}e^{-2\pi b u}-e^{-2\pi b}\right)}{(j-1)!(k-j)!}(1-\cos{2\pi n (1-u)})\cot{\pi(1-u)}\,du \text{,}
\end{equation}
\noindent and therefore, since the infinities now cancel out, one concludes that:
\begin{equation} \label{eq:Hurwitz_im_part}
\Im{\left(\sum_{j=0}^{\infty}\frac{1}{(\ii j+b)^k}\right)}=-\frac{(2\pi)^k}{2}\int_{0}^{1}\sum_{j=1}^{k}\frac{\left(\delta_{1j}+\mathrm{Li}_{-j+1}(e^{-2\pi b})\right)\left(u^{k-j}e^{-2\pi b u}-e^{-2\pi b}\right)}{(j-1)!(k-j)!}\cot{\pi u}\,du 
\end{equation}
\indent The imaginary part does not converge when $k=1$, since the above integral does not converge.

\section{Hurwitz zeta function}
The relation \eqrefe{relacao} implies that when $b$ is a real or purely imaginary number, the created formulae allow to separate the real and imaginary parts of $\zeta(k,-\ii\,b)$.\\

Coincidentally, when the two parts are combined the formula holds even when $b$ is not real. Hence, a simple transformation is applied to yield a more straightforward formulation of the Hurwitz zeta function. For integer $k\ge 2$:
\begin{multline} \label{eq:Hurwitz_form}
\zeta(k,b)=\frac{1}{2b^k}+\frac{(2\pi \ii)^k\left(\delta_{1k}+\mathrm{Li}_{-k+1}(e^{-2\pi\ii\,b})\right)}{4(k-1)!}+\frac{(2\pi\ii)^ke^{-2\pi\ii\,b}}{4}\sum_{j=1}^{k}\frac{\delta_{1j}+\mathrm{Li}_{-j+1}(e^{-2\pi\ii\,b})}{(j-1)!(k-j)!}\\-\frac{\ii(2\pi\ii)^k}{2}\int_{0}^{1}\sum_{j=1}^{k}\frac{\left(\delta_{1j}+\mathrm{Li}_{-j+1}(e^{-2\pi\ii\,b})\right)\left(u^{k-j}e^{-2\pi\ii\,b\,u}-e^{-2\pi\ii\,b}\right)}{(j-1)!(k-j)!}\cot{\pi u}\,du 
\end{multline} 

\subsection{A particular case}
The below closed form was not found in the literature (including the software Mathematica), though it is unlikely that it is really unknown:
\begin{equation} \nonumber
\zeta\left(2,\frac{5}{4}\right)= -16+\pi^2+8G \text{, where } G\text{ is Catalan's constant.}
\end{equation}

\section{The $\zeta(k,b)$ generating function}
In this section, assuming again that $b$ is real, a generating function for $\HP_k(n)$ is derived and its limit when $n$ goes to infinity is determined:
\begin{multline} \nonumber
\sum_{k=1}^{\infty}x^k\sum_{j=1}^{n}\frac{1}{(\ii j+b)^k}=\sum_{k=1}^{\infty}\Big(-\frac{x^k}{2b^k}+\frac{x^k}{2(\ii n+b)^k}+\\ \frac{1}{2}\int_{0}^{1}\sum_{j=1}^{k}\frac{\left(\delta_{1j}+\mathrm{Li}_{-j+1}(e^{-2\pi b})\right)(2\pi x)^j}{(j-1)!}\frac{(2\pi x(1-u))^{k-j}}{(k-j)!}e^{-2\pi b(1-u)}\left(\sin{2\pi n u}+\ii(1-\cos{2\pi n u})\right)\cot{\pi u}\,du \Big)
\end{multline}
\indent The summation within the integral is the general term of the series expansion of the product of the two functions below (per the Leibniz rule):
\begin{equation} \nonumber
p(x)=-\frac{2\pi x}{\left(e^{2\pi x}-e^{2\pi b}\right)}\cdot e^{2\pi x(1-u)}
\end{equation}
\indent Therefore, the power series within the integral can be replaced by the product function, and the formula can be rewritten as:
\begin{multline} \nonumber
\sum_{k=1}^{\infty}x^k\sum_{j=1}^{n}\frac{1}{(\ii j+b)^k}=\sum_{k=1}^{\infty}x^k\left(-\frac{1}{2b^k}+\frac{1}{2(\ii n+b)^k}\right)
\\-\frac{\pi x}{e^{2\pi x}-e^{2\pi b}}\int_{0}^{1}e^{2\pi x(1-u)}e^{2\pi b u}\left(\sin{2\pi n u}+\ii(1-\cos{2\pi n u})\right)\cot{\pi u}\,du 
\end{multline}
\indent Since we want the limit of this expression as $n$ goes to infinity, the harmonic progression (of order 1) must be subtracted from it (the method will show once more why this is the case, subsequently). Thus one has:
\begin{multline} \label{eq:GF_partial}
\sum_{k=2}^{\infty}x^k\sum_{j=1}^{n}\frac{1}{(\ii j+b)^k}=\sum_{k=2}^{\infty}x^k\left(-\frac{1}{2b^k}+\frac{1}{2(\ii n+b)^k}\right)
\\-\pi x\int_{0}^{1}\left(\frac{1}{e^{2\pi b}-1}+\frac{e^{2\pi x(1-u)}}{e^{2\pi x}-e^{2\pi b}}\right)e^{2\pi b u}\left(\sin{2\pi n u}+\ii(1-\cos{2\pi n u})\right)\cot{\pi u}\,du 
\end{multline}

\subsection{The real part}
Let us consider the limit of the real part first. By following the same thought process that has been laid out in section \secrefe{Real}, one concludes that the real part of the limit of \eqrefe{GF_partial} is given by:
\begin{equation} \label{eq:GF_real_part}
\Re{\left(\sum_{k=2}^{\infty}x^k\sum_{j=1}^{\infty}\frac{1}{(\ii j+b)^k}\right)}=\frac{x^2}{2b(x-b)}+\frac{\pi x(e^{2\pi x}-1)}{(e^{-2\pi b}-1)(e^{2\pi x}-e^{2\pi b})}
\end{equation}

\subsection{The imaginary part}
Note that when this demonstration was first devised, we still had not had the insight that was used to derive the limit of the imaginary part from section \secrefe{approx}, so the following rationale may be unnecessarily convoluted.\\

Regarding the imaginary part of the limit of \eqrefe{GF_partial}, the two terms outside of the integral can obviously be discarded, as one is real and the other one goes to zero as $n$ goes to infinity, leaving one with:
\begin{multline} \nonumber
\Im{\left(\sum_{k=2}^{\infty}x^k\sum_{j=1}^{\infty}\frac{1}{(\ii j+b)^k}\right)}=\lim_{n\to\infty}-\pi x\int_{0}^{1}\left(\frac{1}{e^{2\pi b}-1}+\frac{e^{2\pi x(1-u)}}{e^{2\pi x}-e^{2\pi b}}\right)e^{2\pi b u}(1-\cos{2\pi n u})\cot{\pi u}\,du 
\end{multline}
\indent Now one must find the limit of the integral as $n$ approaches infinity. To solve the two limits at once, let us use $c=c(x)$ as the coefficient of $u$:
\begin{equation} \nonumber
\lim_{n\to\infty}\int_{0}^{1}e^{-c\,u}(1-\cos{2\pi n u})\cot{\pi u}\,du 
\end{equation}
\indent To ensure that nothing is missed, let us make our initial formula more aligned with the asymptotic formulae from section \secrefe{approx} by changing $u$ for $1-u$:
\begin{multline} \nonumber
\Im{\left(\sum_{k=2}^{\infty}x^k\sum_{j=1}^{n}\frac{1}{(\ii j+b)^k}\right)}\\=\lim_{n\to\infty}-\pi xe^{2\pi b}\int_{0}^{1}\left(\frac{e^{-2\pi bu}}{e^{2\pi b}-1}+\frac{e^{2\pi(x-b)u}}{e^{2\pi x}-e^{2\pi b}}\right)(1-\cos{2\pi n(1-u)})\cot{\pi(1-u)}\,du 
\end{multline}
\indent Also, a variable $y$ is introduced into the integral that needs to evaluated, to help with the power series manipulations that need to be performed. In the end, one needs to set $y$ to 1 (and $c$ to $2\pi b$ or $2\pi(x-b)$):
\begin{multline} \nonumber
\int_{0}^{1}e^{-c y\,u}(1-\cos{2\pi n(1-u)})\cot{\pi(1-u)}\,du=
\\ \int_{0}^{1}\left(1+\sum_{k=0}^{\infty}\frac{(-c y)^{2k+1}}{{(2k+1)}!}u^{2k+1}+\frac{(-c y)^{2k+2}}{{(2k+2)}!}u^{2k+2}\right)(1-\cos{2\pi n(1-u)})\cot{\pi(1-u)}\,du
\end{multline}
\indent Note the constant 1 can be ignored, since that integral is zero for all integer $n$. After each $u^k(1-\cos{2\pi n(1-u)})\cot{\pi(1-u)}\,du$ is replaced with their asymptotic formulae from section \secrefe{approx}, part of the above expression reduces to:
\begin{equation}
\frac{\gamma+\log{n}}{\pi}\sum_{k=0}^{\infty}\left(-\frac{c^{2k+1}}{(2k+1)!}+\frac{c^{2k+2}}{(2k+2)!}\right)=\frac{\gamma+\log{n}}{\pi}\left(-\sinh{c}-1+\cosh{c}\right)\nonumber
\end{equation}
\indent This part diverges, and it is due to $\HP(n)$, as mentioned at the beginning, which was subtracted from the generating function. Therefore, these infinities are expected to cancel out when the terms that have $c=2\pi b$ and $c=-2\pi(x-b)$ are added to our initial formula.\\

Now, the below is the part that does not diverge:
\begin{equation} \label{eq:series}
\frac{1}{\pi}\sum_{k=0}^{\infty}\left(-(c y)^{2k+1}\sum_{j=1}^{k}\frac{(-1)^j(2\pi)^{-2j}}{(2k+1-2j)!}\zeta{(2j+1)}+(c y)^{2k+2}\sum_{j=1}^{k}\frac{(-1)^j(2\pi)^{-2j}}{(2k+2-2j)!}\zeta{(2j+1)}\right)
\end{equation}
\indent It is not a very simple series, and it is not straightforward to know how to proceed from here. Herein, the integral representation that was derived for $\zeta(2j+1)$ in \citena{GHNR} can be used. Below it has been slightly modified:
\begin{equation} \nonumber
\zeta(2j+1)=-\frac{(-1)^j(2\pi)^{2j+1}}{2}\int_{0}^{1}\sum_{q=0}^{j}\frac{B_{2 q}\left(2-2^{2 q}\right)u^{2j-2 q+1}}{(2 q)!(2j-2 q+1)!}\cot{\pi u}\,du 
\end{equation}
\indent For simplicity, let us start with the first part of \eqrefe{series}:
\begin{equation} \nonumber
-\frac{1}{\pi}\sum_{k=0}^{\infty}(c y)^{2k+1}\sum_{j=1}^{k}\frac{(-1)^j(2\pi)^{-2j}}{(2k+1-2j)!}\left(-\frac{(-1)^j(2\pi)^{2j+1}}{2}\int_{0}^{1}\sum_{q=0}^{j}\frac{B_{2q}\left(2-2^{2q}\right)u^{2j-2q+1}}{(2q)!(2j-2q+1)!}\cot{\pi u}\,du\right) \Rightarrow
\end{equation}
\begin{multline} \nonumber
\sum_{k=0}^{\infty}(c y)^{2k+1}\sum_{j=1}^{k}\frac{1}{(2k+1-2j)!}\int_{0}^{1}\sum_{q=0}^{j}\frac{B_{2q}\left(2-2^{2q}\right)u^{2j-2q+1}}{(2q)!(2j-2q+1)!}\cot{\pi u}\,du=
\\ \frac{1}{cy}\int_{0}^{1}\sum_{k=0}^{\infty}\left(\sum_{j=0}^{k}\frac{(c y)^{2k+1-2j}}{(2k+1-2j)!}(c y)^{2j+1}\sum_{q=0}^{j}\frac{B_{2q}\left(2-2^{2q}\right)u^{2j-2q+1}}{(2q)!(2j-2q+1)!}-\frac{(c y)^{2k+2}u}{(2k+1)!}\right)\cot{\pi u}\,du 
\end{multline}
\indent Again, the above power series is the product of two slightly familiar functions. One of them is obviously $\sinh{cy}$, and the other one (which appeared in \citena{GHNR} in non-hyperbolic form) is:
\begin{equation} \nonumber
\sum_{j=0}^{\infty}(c y)^{2j+1}\sum_{q=0}^{j}\frac{B_{2q}\left(2-2^{2q}\right)u^{2j-2q+1}}{(2q)!(2j-2q+1)!}=cy\csch{cy}\sinh{cuy}
\end{equation}
\indent Therefore, the final conclusion is:
\begin{equation} \nonumber
-\frac{1}{\pi}\sum_{k=0}^{\infty}(c y)^{2k+1}\sum_{j=1}^{k}\frac{(-1)^j(2\pi)^{-2j}}{(2k+1-2j)!}\zeta{(2j+1)}=\sinh{cy}\int_{0}^{1}\left(\csch{cy}\sinh{cuy}-u\right)\cot{\pi u}\,du
\end{equation}
\indent And since the second part follows an analogous thought process, its development is omitted, but the end result is below:
\begin{equation} \nonumber
\frac{1}{\pi}\sum_{k=0}^{\infty}(c y)^{2k+2}\sum_{j=1}^{k}\frac{(-1)^j(2\pi)^{-2j}}{(2k+2-2j)!}\zeta{(2j+1)}=(1-\cosh{cy})\int_{0}^{1}\left(\csch{cy}\sinh{cuy}-u\right)\cot{\pi u}\,du
\end{equation}
\indent Finally, one can sum it all up by making $y=1$:
\begin{multline} \nonumber
\frac{1}{\pi}\sum_{k=0}^{\infty}-c^{2k+1}\sum_{j=1}^{k}\frac{(-1)^j(2\pi)^{-2j}}{(2k+1-2j)!}\zeta{(2j+1)}+c^{2k+2}\sum_{j=1}^{k}\frac{(-1)^j(2\pi)^{-2j}}{(2k+2-2j)!}\zeta{(2j+1)}=\\-(e^{-c}-1)\int_{0}^{1}\left(\csch{c}\sinh{cu}-u\right)\cot{\pi u}\,du
\end{multline}
\indent To conclude, one must evaluate the initial expression with the above identity. After that is carried out, one finds that:
\begin{multline} \label{eq:GF_im_part}
\Im{\left(\sum_{k=2}^{\infty}x^k\sum_{j=1}^{\infty}\frac{1}{(\ii j+b)^k}\right)}\\=\pi x\int_{0}^{1}\left(\csch{2\pi(x-b)}\sinh{2\pi(x-b)u}-\sinh{2\pi b}\csch{2\pi bu}\right)\cot{\pi u}\,du
\end{multline}

\subsection{Conclusion}
\indent If $b$ is not assumed to be real then, provided that $\ii\,b\notin \mathbb{Z}$ (and provided that the right-hand side does not contain singularities, such as $b=x$), one can simply state that, for $b,x\in\mathbb{C}$, when the limit exists it is given by:
\begin{multline} \nonumber
\sum_{k=2}^{\infty}x^k\sum_{j=1}^{\infty}\frac{1}{(\ii j+b)^k}=\frac{x^2}{2b(x-b)}+\frac{\pi x(e^{2\pi x}-1)}{(e^{-2\pi b}-1)(e^{2\pi x}-e^{2\pi b})}+\\ \ii \pi x\int_{0}^{1}\left(\csch{2\pi(x-b)}\sinh{2\pi(x-b)u}-\csch{2\pi b}\sinh{2\pi bu}\right)\cot{\pi u}\,du
\end{multline}
\indent However, when the left-hand side diverges, the expression on the right may still converge, if it does not have singularities --- which means that it is an analytic continuation of the left-hand side.\\

The above can be turned into a better looking equation without non-real numbers, which in principle holds if $b\neq 0$, $x \neq b$ and $2b$ and $2(x-b)$ are not integers:
\begin{multline} \label{eq:GF_final}
f(x)=\sum _{k=2}^{\infty}x^k\left(\zeta(k,b)-\frac{1}{b^k}\right)=\frac{x^2}{2b(x-b)}-\frac{\pi x\sin{\pi x}}{2\sin{\pi b}}\csc{\pi(x-b)}\\-\pi x\int _0^1\left(\frac{\sin{2\pi(x-b)u}}{\sin{2\pi(x-b)}}-\frac{\sin{2\pi b u}}{\sin{2\pi b}}\right)\cot{\pi u}\,du
\end{multline}
\indent Note that $x=2b$ causes the integral to vanish, which means that $f(x)$ has a closed form at that point.\\

To know what the formula becomes when $b$ is a positive integer, one can rely on the below identity:
\begin{equation} \nonumber
f(x)=\sum _{k=2}^{\infty}x^k\sum _{j=1}^{\infty}\frac{1}{(j+b)^k}=\sum _{k=2}^{\infty}x^k(\zeta{(k)}-H_{k}(b))
\end{equation}
\indent In general, when $b$ is an integer, $f(x)$ can be obtained by using the generating functions created for $H_k(n)$ and $\zeta(k)$ in \citena{GHNR}, which, after all the necessary calculations are performed, leads to:
\begin{equation} \nonumber
f(x)=\begin{cases} 
\frac{1}{2}-\frac{\pi x}{2}\cot{\pi x}-\pi x\int _0^1\left(\frac{\sin{2 \pi xu}}{\sin{2\pi x}}-u\right)\cot{\pi u}\,du & \text{, if }b=0 \\
\frac{x^2}{2b(x-b)}-\frac{\pi x}{2}\cot{\pi x}-\pi x\int _0^1\left(\frac{\sin{2\pi (x-b)u}}{\sin{2\pi x}}-u\cos{2\pi b u}\right)\cot{\pi u}\,du & \text{, if }b \in\mathbb{Z_+}\\
1+\frac{x^2}{2b(x-b)}-\frac{\pi x}{2}\cot{\pi x}-\pi x\int _0^1\left(\frac{\sin{2\pi(x-b)u}}{\sin{2\pi x}}-u\cos{2\pi b u}\right)\cot{\pi u}\,du & \text{, if }b \in\mathbb{Z_-} \text{,}
\end{cases} 
\end{equation}
\noindent where singularities are being disregarded (in the case of negative integer $b$), which leaves only the half-integers $b$ unaccounted for.\\

Note some of these functions have a removable singularity at zero, so they are analytic at zero. So, for $k\geq 2$, $\zeta(k,b)$ can also be obtained as:
\begin{equation}\nonumber
\zeta(k,b)=\frac{1}{b^k}+\frac{f^{(k)}(0)}{k!}
\end{equation}



\begin{thebibliography}{1}

\bibitem{Abramo} M. Abramowitz, I. A. Stegun, {\em Handbook of Mathematical Functions with Formulas, Graphs, and Mathematical Tables (9th printing ed.), New York: Dover,}  1972.

\bibitem{GHNR} Risomar Sousa, Jose {\em Generalized Harmonic Numbers, eprint arXiv:1810.07877,} 2018.

\bibitem{GHP} Risomar Sousa, Jose {\em Generalized Harmonic Progression, eprint arXiv:1811.11305,} 2018.

\bibitem{CHP} Risomar Sousa, Jose {\em Generalized Harmonic Progression Part II, eprint arXiv:1902.01008,} 2019.
\end{thebibliography}
\end{document}